\documentclass[final,leqno]{siamltex}
\usepackage{amsmath}
\usepackage{amsfonts}
\usepackage{color}
\usepackage[latin1]{inputenc}
\usepackage{algpseudocode}
\usepackage{algorithm}
\usepackage{multirow}
\usepackage{caption}
\usepackage{lscape}
\usepackage[update,prepend]{epstopdf} 
\usepackage{graphicx}
\usepackage{multirow}
\usepackage{url}
\usepackage{diagbox}

\newcommand{\be}{\begin{eqnarray}}
\newcommand{\e}{\end{eqnarray}}
\newcommand{\bes}{\begin{eqnarray*}}
\newcommand{\es}{\end{eqnarray*}}
\newcommand{\beq}{\begin{equation}}
\newcommand{\eeq}{\end{equation}}

\newcommand{\trace}[1]{\ensuremath{\mathop{\mathrm{trace}}\left( #1 \right)}}
\renewcommand{\vec}[1]{\ensuremath{\mathop{\mathrm{vec}}\left( #1 \right)}}

\author{Peter Benner, Akwum Onwunta \thanks{Corresponding author}, and Martin Stoll\\
Max Planck Institute for Dynamics of Complex Technical Systems, Sandtorstrasse 1, 39106 Magdeburg, Germany }

\begin{document}
\title{An Inexact Newton-Krylov method for stochastic eigenvalue problems}
\author{Peter Benner\footnotemark[1], Akwum Onwunta\footnotemark[2] \and Martin Stoll\footnotemark[3]}
\renewcommand{\thefootnote}{\fnsymbol{footnote}}
\footnotetext[1]{Computational Methods in Systems and Control Theory, Max Planck Institute
for Dynamics of Complex Technical Systems,
Sandtorstrasse 1,
39106 Magdeburg,
Germany,
(\url{benner@mpi-magdeburg.mpg.de})}
\footnotetext[2]{Corresponding author; Computational Methods in Systems and Control Theory, Max Planck Institute
for Dynamics of Complex Technical Systems,
Sandtorstrasse 1,
39106 Magdeburg,
Germany,
(\url{onwunta@mpi-magdeburg.mpg.de})}
\footnotetext[3]{Numerical
Linear Algebra for Dynamical Systems Group,
Max Planck Institute for Dynamics of Complex Technical Systems,
Sandtorstrasse 1,
39106 Magdeburg,
Germany,
(\url{stollm@mpi-magdeburg.mpg.de});
Technische Universit\"{a}t Chemnitz, Faculty of Mathematics, Professorship Scientific Computing, 09107 Chemnitz, Germany,
(\url{martin.stoll@mathematik.tu-chemnitz.de})}
\renewcommand{\thefootnote}{\arabic{footnote}}
\maketitle

\begin{abstract}
This paper aims at the efficient numerical solution of stochastic eigenvalue problems. 
Such problems often lead to prohibitively
high dimensional  systems with tensor product structure  when discretized with the stochastic 
 Galerkin method. Here, we exploit this inherent tensor product structure to develop a
globalized low-rank  inexact Newton method with which we tackle the stochastic eigenproblem. 
We illustrate  the effectiveness of our solver with numerical experiments. 
\end{abstract}

\begin{keywords}
Stochastic Galerkin system, Krylov methods, eigenvalues, eigenvectors, low-rank solution,  preconditioning. \\
\end{keywords}

\begin{AMS}35R60, 60H15, 60H35, 65N22, 65F10, 65F50\end{AMS}

\pagestyle{myheadings}
\thispagestyle{plain}
\markboth{Efficient solvers for stochastic eigenvalue problems}{Efficient solvers for stochastic eigenvalue problems}
\section{Introduction}
\label{sec1}
In many areas of computational science and engineering, eigenvalue problems play an 
important role. This is, for example, the case in structural
mechanics, where eigenvalue problems typically appear in the context of vibrations and buckling.
For deterministic problems, there are currently well-established   algorithms 
dedicated to the computation of eigenvalues and eigenvectors, see, e.g., \cite{GL96}.  However,
in many cases of practical  interest, physical characteristics are not always completely deterministic.
For instance, the stiffness of a plate can locally be reduced by material imperfections, or the
velocity of a flow can be influenced by turbulence. In recent times, an increasingly important way to model such problems is
by describing the uncertain problem characteristics more
realistically using random variables. By  doing  so, one would then gain  more insight 
regarding the effect of the uncertainties on the model. 
This approach then leads  to a stochastic eigenvalue problem (SEVP).


It is worth pointing out that the consequence of modeling  
the input parameters of a physical problem as random variables is
that the desired output  naturally inherits the stochasticity in the model.
Generally speaking, there are  two  broad techniques for analyzing and quantifying uncertainty in a given model:
simulation-based methods and expansion-based methods.
In the simulation- (or sampling-) based methods, 
the stochastic moments of the eigenvalues
and eigenvectors are obtained by generating ensembles of random realizations for the prescribed
random inputs and utilizing repetitive deterministic solvers for each realization.
Prominent among this class of methods is the classical Monte Carlo method.
This method has been applied to many problems and its implementations
are straightforward. It is (formally) independent of the dimensionality
of the random space; that is, it is independent of the number of random variables used to
characterize the random inputs.
It does, however,  exhibit a very slow  convergence rate \cite{TKXP12}.  To accelerate its convergence, several techniques
have been developed: the multilevel Monte Carlo method \cite{CGST11}, the
quasi-Monte Carlo  method \cite{N92},  the Markov chain Monte Carlo
method  \cite{GRS95},  etc.
Although these methods  can improve the efficiency of the traditional
Monte Carlo method, additional restrictions are imposed based on
their specific designs and their applicability is limited.

The expansion-based  methods for uncertainty 
analysis and quantification are often designed to retain the advantages of Monte Carlo
simulations; in particular, they enable one to compute the full statistical characteristics of the
solution, while reducing the simulation time. A typical example
of the expansion-based methods are the  spectral  stochastic finite element methods (SFEM) \cite{GS91, PE09}; they   rely on the approximation of  the random eigenvalues and eigenvectors
by projecting them onto a global basis and are considerably less expensive than the simulation-based methods. 
We will, in particular, employ mainly SFEM in this paper.


During the last two decades, there has been a lot of research on 
SFEM for uncertainty analysis and quantification for solutions of partial differential equations \cite{BOS14, BOS15a, PE09}.
However, SFEM 
  for SEVPs has been so far   much less addressed in the literature.
 To a great extent, most research on SEVPs has, in fact, focused more
on  simulation-based techniques \cite{PSS02, SS01}.
Nevertheless, relatively few attempts have been made to approximate the stochastic moments of both the eigenvalues
and eigenvectors through the use of spectral methods \cite{GG07, HKL15, VGH06}. 
In \cite{VGH06}, the authors propose  algorithms based on the inverse power method together with spectral methods 
for computing approximate eigenpairs of both symmetric and non-symmetric SEVPs.
The method proposed in \cite{GG07} essentially
 rewrites the eigenvalue problem resulting from a spectral discretization
 (which we henceforth refer to as stochastic Galerkin method (SGM)) as a set of nonlinear equations with
 tensor product structure, 
 which are then solved 
 using the Newton-Raphson method.  In the spirit of \cite{GG07}, this paper presents an algorithm to determine the spectral
expansions of the eigenvalues and the eigenvectors based on a Newton's method and SGM. However, unlike
\cite{GG07}, this work specifically focuses on the use of a {\it globalized low-rank   inexact Newton method}
to tackle the eigenproblem.  

Now, recall that under certain conditions, the iterates produced by the Newton's
method converge quadratically to a solution $x^{\ast}$ of a given nonlinear system, and  those of the inexact
Newton method can obtain super-linear convergence \cite{AML2007, EW96, STW97}. Both cases, however,
assume an initial guess $x_0$ sufficiently close to $x^{\ast}$.  Generally speaking, globalizing the
inexact Newton method  means augmenting the method with additional conditions on the choices
of iterates $\{x_k\}$ to enhance the likelihood of convergence to $x^{\ast}$, 
see e.g. \cite{STW97} for details of different globalization
techniques. 
The advantages of globalization notwithstanding\footnote{It is important to note that no  globalization strategy  determines a sequence that
converges to a solution for every problem; rather,  globalization techniques are essentially used
only to enhance the likelihood of convergence to some solution of the problem.}, 
a drawback of Newton-type methods is that for  fairly large eigenproblems
such as the SEVPs considered in this work, they require
considerable computational effort to solve the linear system arising from each Newton step.
The aim of this paper is therefore to mitigate this computational challenge by 
exploiting the inherent tensor product structure in the SEVP to  tackle the stochastic eigenproblem.
More precisely, we combine low-rank Krylov solvers with a globalized inexact Newton method
to efficiently solve  SEVPs.

The rest of the paper is organized as follows. In Section \ref{probst}, we present the problem that
we would like to solve in this paper.  Next, Section \ref{sgm} gives an overview of the stochastic
Galerkin method on which we shall rely to discretize our model problem. After discussing our 
globalized low-rank  inexact Newton solver in Section \ref{newtonkry},  we proceed 
to Section \ref{numresults} to provide
the numerical results to buttress the efficiency of the proposed solver, while Section \ref{conc} 
draws some conclusions on the findings in this work.

\section{Problem statement}
\label{probst}
 Let the triplet $(\Omega,\mathcal{F},\mathbb{P})$ denote
a complete probability space,  where $\Omega$
is  the set of elementary events, $\mathcal{F}\subset 2^{\Omega}$ is a $\sigma$-algebra on $\Omega$
 and $\mathbb{P}:\mathcal{F}\rightarrow [0,1]$ is an appropriate
probability measure. Let $\mathcal{D}\subset \mathbb{R}^d$ with $d\in\{1,2,3\},$ be a bounded physical domain.
In this paper, we consider  the following eigenvalue problem for an $N_x$-dimensional real symmetric random
matrix 
\be
\label{randeqn}
\mathcal{A}(\omega) \varphi(\omega) = \lambda(\omega) \varphi(\omega),
\e
subject to the normalization condition
\be
\label{normeig}
\varphi(\omega)^T\varphi(\omega) =1, 
\e
where 
\[
 \lambda(\omega)\in \mathbb{R}, \quad \varphi(\omega)\in \mathbb{R}^{N_x}, \quad
 \mathcal{A}(\omega)\in \mathbb{R}^{{N_x}\times {N_x}}, \quad \omega \in \Omega.
\]
The matrix $\mathcal{A}(\omega)$ represents, for example, the stiffness matrix in a structural mechanics
problem \cite{GG07}.
In this case, the stochasticity in $\mathcal{A}(\omega)$  is often   inherited from the randomness in the underlying
physical system such as elastic and dynamic parameters. Moreover, we assume that the randomness in the model is induced 
by a prescribed finite number of random variables $\xi:=\{\xi_1,\xi_2,\ldots, \xi_m\},$ where
$ m\in\mathbb{N}$ and $\xi_i(\omega):{\Omega}\rightarrow \Gamma_i\subseteq \mathbb{R}.$
We also make the simplifying assumption that
each random variable is independent and characterized by a probability density function $\rho_i:\Gamma_i\rightarrow[0,1].$
If the distribution measure of the random vector $\xi(\omega)$ is absolutely continuous with respect
to the Lebesgue measure, then there exists a joint probability density function
$\rho:\Gamma\rightarrow \mathbb{R}^+,$ where
$\rho(\xi)=\prod^m_{i=1}\rho_i(\xi_i), $
and $\rho\in L^{\infty}(\Gamma).$
Furthermore, we can now replace the probability space $(\Omega,\mathcal{F},{\mathbb{P}})$ 
with $(\Omega,\mathbb{B}(\Gamma),\rho(\xi)d\xi),$  where $\mathbb{B}(\Gamma)$ denotes the
Borel $\sigma$-algebra on $\Gamma$ and $\rho(\xi)d\xi$ is the finite measure of the vector $\xi.$
Then, the expected value of the product of measurable functions on $\Gamma$ determines the Hilbert space
${L}_\rho^2(\Omega,\mathbb{B}(\Gamma),\rho(\xi)d\xi),$ with inner product
\[
 \left<u,v\right>:=\mathbb{E}[uv]=\int_\Gamma u(\xi)v(\xi)\rho(\xi)d\xi,
\]
where the symbol $\mathbb{E}$ denotes mathematical expectation.  

In this paper, we assume that the random matrix $\mathcal{A}(\omega)$ in (\ref{randeqn}) admits the representation
 \be
 \label{rmatrix}
 \mathcal{A}(\omega) = A_0 + \sum_{k=1}^m \xi_k(\omega) A_k, \quad m\in\mathbb{N}, \;\;  A_k\in \mathbb{R}^{{N_x}\times {N_x}},
  \; \; k =0,1,\ldots, m,
 \e
where $\{\xi_k\}$ are independent  random variables. This is indeed the case if a 
Karhunen-Lo\`{e}ve  expansion (KLE) is used to discretize random stiffness properties; see, e.g., \cite{HKL15, PE09, AO16}.
Furthermore, the stochastic eigenvalues and eigenvectors in this work are
approximated using  the so-called
{\it generalized polynomial chaos expansion} (gPCE)
\cite{BOS14, AO16, XS09}. More precisely, the $\ell$th random eigenvalue
and eigenvector are given, respectively, as 
\be
\label{lambdas}
\lambda_\ell(\omega) = \sum_{k=0}^{N_{\xi}-1} \lambda^{(\ell)}_k \psi_k(\xi(\omega)), \quad \lambda^{(\ell)}_k\in \mathbb{R},
\e
and 
\be
\label{phis}
\varphi_\ell(\omega) = \sum_{k=0}^{N_{\xi}-1} \varphi^{(\ell)}_k \psi_k(\xi(\omega)), \quad \varphi^{(\ell)}_k\in \mathbb{R}^{N_x},
\e
where $\{\psi_i\}$ are multidimensional Legendre basis polynomials expressed as functions of the random vector $\xi,$ with properties
\[
 \mathbb{E}(\psi_k) =\delta_{k0}  \quad \; \mbox{and} \quad \; 
 \mathbb{E}(\psi_j\psi_k) =\delta_{jk}\mathbb{E}(\psi_k^2).
\]
The  spectral expansions (\ref{lambdas}) and (\ref{phis}) are the   gPCE 
of the  random quantities $\lambda_\ell(\omega)$ and  $\varphi_\ell(\omega),$ respectively.
Throughout this paper, we use normalized Legendre basis polynomials in which case $\mathbb{E}(\psi_i^2) =1,$
so that $\mathbb{E}(\psi_i\psi_j) =\delta_{ij}.$
We remark here that $N_{\xi}$ in (\ref{lambdas})  and  (\ref{phis}) is chosen in such a way that $N_{\xi} > m.$ In particular, using total degree Legendre
polynomials $\psi_i$ yields 
\be
\label{Nxi}
N_{\xi}=(m+r)!/m!r!,
\e
where $r$ is the degree of $\psi_i,$ see e.g. \cite{PE09}.

In what follows,  we will, for notational convenience,  omit the index $\ell$ associated with the $\ell$th eigenpair. 
It is  pertinent to note here the difference between the structure of a deterministic and a random eigenproblem.
In the deterministic case, a typical eigenpair is of the form $(\lambda,\varphi),$ where 
$\lambda\in \mathbb{R}$ and $\varphi\in \mathbb{R}^{N_x},$ with ${N_x}$ denoting the size of the 
deterministic matrix $\mathcal{A}$. 
In the stochastic case, however, the eigenpair corresponding to $\ell$th physical mode consists of the set
\be
\label{eigenpair}
x:= \{\lambda_0,\lambda_1,\ldots, \lambda_{N_{\xi}-1},
 \varphi_0, \varphi_1, \ldots, \varphi_{N_{\xi}-1}\}.
\e


\section{Stochastic Galerkin method} 
\label{sgm}
The stochastic Galerkin method is based on the projection
\be
\label{projn}
 \left<\mathcal{A}\varphi, \psi_k\right> =  \left<\lambda \varphi, \psi_k\right>, \quad k = 0, \ldots, N_{\xi}-1,
 \quad \ell=1,\ldots {N_x}.
\e
Substituting (\ref{rmatrix}), (\ref{lambdas}), and (\ref{phis}) into (\ref{projn}) yields the nonlinear algebraic 
equations
\be
\sum_{i=0}^{m-1}\sum_{j=0}^{N_{\xi}-1} \mathbb{E}(\xi_i\psi_j\psi_k) A_i \varphi_j =
\sum_{i=0}^{N_{\xi}-1}\sum_{j=0}^{N_{\xi}-1} \mathbb{E}(\psi_i\psi_j\psi_k) \lambda_i \varphi_j,
\; k = 0, \ldots, N_{\xi}-1,
\e
which can be rewritten in Kronecker product notation as
\be
\label{matrices}
\underbrace{\left[G_0 \otimes A_0 +\sum_{k=1}^{m} G_k \otimes A_k\right]}_{:=A} \Phi = 
\left[\sum_{k=0}^{N_{\xi}-1}\lambda_k \underbrace{(H_k \otimes  {\bf I})}_{:=B_k} \right] \Phi,
\e
where ${\bf I}$ is the  identity matrix and
\be
 \label{Gs}
 \begin{cases}
   G_0= \mbox{diag}\left( \left\langle \psi^2_0  \right \rangle, \left\langle \psi^2_1 \right \rangle,\ldots,  \left\langle \psi^2_{N_{\xi}-1}  \right \rangle \right),\\
   G_k(i,j) = \left< \psi_i\psi_j\xi_k\right>,\;\;  k=1,\ldots, m,\\
   H_k(i,j) = \left< \psi_i\psi_j\psi_k\right>,\;\; k=0,\ldots, N_{\xi}-1,\\
   \Phi = \left(\varphi_0, \varphi_1, \ldots, \varphi_{N_{\xi}-1} \right)\in \mathbb{R}^{{N_x}N_{\xi}}.
 \end{cases}
\e
Here, the block $A_0$ (as well as  $A$ itself) is symmetric and positive definite; it captures the mean information in the model
and appears on the diagonal blocks of $A,$ whereas
the other blocks $A_k, \;k=1,\ldots, m,$ represent the fluctuations in the model. 
Moreover,  the random variables $\{\xi_k\}^m_{k=1}$ are centered, normalized and  independent; see e.g., \cite{PE09}.

Recalling that $N_{\xi}> m,$ we see that (\ref{matrices}) can also be expressed as
\be
\label{newmat}
\underbrace{\sum_{k=0}^{N_{\xi}-1} \left[  (G_k \otimes A_k) -
\lambda_k (H_k \otimes  {\bf I}) 
\right] \Phi}_{:=E} = 0, \quad  G_k = A_k= 0,\; \mbox{for}\;  k>m.
\e

Now, observe that the problem (\ref{matrices}) can be considered as an {\it eigentuple-eigenvector} problem:
\be
\label{eigetuple}
 A\Phi = \sum_{k=0}^{N_{\xi}-1} \lambda_k B_k \Phi,
\e
in which one needs to find an eigentuple $\Lambda:=(\lambda_0, \ldots, \lambda_{N_{\xi}-1})\in \mathbb{R}^{N_{\xi}}$ and an 
eigenvector $\Phi  \in   \mathbb{R}^{{N_x}N_{\xi}},$ 
where $A:=\sum_{k=0}^{m} G_k \otimes A_k$ and $B_k:=\lambda_k (H_k \otimes  {\bf I}).$
Note that $B_0:=H_0 = G_0= {\bf I}.$ Thus, the case $k=0$ in (\ref{eigetuple}) corresponds to the standard deterministic 
eigenproblem 
\be
\label{onepara}
 A\Phi =  \lambda_0  \Phi,
\e
which has already been studied extensively \cite{Saad2011}.
For $k=1$ (that is, $N_{\xi}=2$), we obtain
\be
\label{twopara}
  (A - \lambda_1B_1)\Phi =  \lambda_0 B_0\Phi,
\e
which yields a standard eigenproblem for each fixed value of $\lambda_1$. Moreover, since $A, B_0$ and $B_1$ are symmetric
matrices (with $B_0$ being positive definite),  we have  a continuum of real solutions
$\lambda_0(\lambda_1)$ parameterized by $\lambda_1.$ The existence of the continuum of real solutions is not surprising
since there are $2{N_x} + 2= 2({N_x}+1)$ unknowns (that is, $\lambda_0, \lambda_1$ and the components of $\Phi$) in only $2{N_x}$ equations.
To circumvent this situation, it is proposed in \cite{GG07}  to prescribe an additional condition via the normalization
of the eigenvectors as in (\ref{randeqn}). This could then make it feasible to determine $\lambda_1$ and thereby reduce the
two-parameter problem (\ref{twopara}) to a one-parameter eigenproblem (\ref{onepara}). Thus,  the existence of a continuum
of real solutions  could make (\ref{twopara}) numerically tractable by permitting its reduction to a sequence of solutions
of (\ref{onepara}), see e.g. \cite{BC78} for details.

The two-parameter eigenproblem has been considered by Hochstenbach  and his co-authors in \cite{HKP05, HP02} following 
a Jacobi-Davidson approach. However, unlike the  problem under consideration in 
this work, for which the resulting system is
coupled, these authors focused on decoupled systems. Moreover, the approach that the authors adopted is quite complicated 
for two-parameter problems and can hardly be applied to multi-parameter eigenproblems considered in this 
paper.  The approach considered here follows closely the framework of 
\cite{GG07}. More specifically,  our method relies on a Newton-Krylov solution technique, which 
we proceed  to discuss 
in Section \ref{newtonkry}.

\newpage
\section{Newton-Krylov approaches}
\label{newtonkry}
\subsection{The Newton system for stochastic eigenvalue problem}
As we already pointed out in Section \ref{sgm}, the problem (\ref{eigetuple}) contains more unknowns than equations.
As suggested  in \cite{GG07}, we incorporate the normalization condition of the eigenvectors
so that the random eigenproblem is posed as  a set of 
\be
\label{Jsize}
{N_x}N_{\xi}+N_{\xi} =({N_x}+1)N_{\xi} 
\e
non-linear deterministic equations for each 
physical mode of the stochastic system. 
To this end, observe that SGM discretization of (\ref{normeig}) yields \cite{GG07}
\be
\label{normalized}
\sum_{i=0}^{N_{\xi}-1} \sum_{j=0}^{N_{\xi}-1} \mathbb{E}(\psi_i\psi_j\psi_k) \varphi^T_i \varphi_j 
= \delta_{k0},\quad k = 0, \ldots, N_{\xi}-1,
\e
or, equivalently,
\be
\label{newnorm}
\Phi^T(H_k \otimes {\bf I})\Phi = \delta_{k0}, \quad k=0,1,\ldots, N_{\xi}-1.
\e

The Newton's method is a well-established iterative method.
For a well-chosen initial iterate, the method exhibits local quadratic convergence.
In this method, (\ref{newmat}) and (\ref{newnorm}) are simultaneously expressed in the form
$F(x)=0,$ where $x= (\Lambda,\Phi)\in \mathbb{R}^{({N_x}+1)N_{\xi}}$ is a vector containing the solution set defined in (\ref{eigenpair}). 
More precisely,  we have
\be
F(x)=\left[\begin{array}{c}
\sum_{k=0}^{N_{\xi}-1} \left[  (G_k \otimes A_k) -
\lambda_k (H_k \otimes  {\bf I}) 
\right] \Phi \\
\Phi^T(H_0 \otimes {\bf I})\Phi -1 \\
\Phi^T(H_1 \otimes {\bf I})\Phi \\
\vdots \\
\Phi^T(H_{N_{\xi}-1} \otimes {\bf I})\Phi \\
\end{array}\right] =
\left[\begin{array}{c}
0 \\
0 \\
\vdots \\
0\\
\end{array}\right].
\e

The Newton iteration for $F(x)=0$ results from a multivariate Taylor expansion about a current point $x_k$ :
\[
F(x_{k+1}) = F(x_k) + F'(x_k)(x_{k+1} - x_k) + \mbox{higher-oder terms}. 
\]
Setting the left-hand side to zero and neglecting the terms of higher-order curvature yields a  Newton
method; that is, given an initial iterate $x_0$, we obtain an iteration over a sequence 
of linear systems (or the Newton equations)
\be
\label{NEq}
F(x_k) + F'(x_k) s_k=0, 
\e
where  $x_k$ is the current iterate. Moreover, $F(x) $  is the 
vector-valued function of nonlinear residuals and  $\mathcal{J}:=F' $ is the associated Jacobian 
matrix, $x$ is the state vector to be found, and $k$ is the  iteration index.
Forming each element of $\mathcal{J}$ requires taking analytic or discrete derivatives of the system of equations with
respect to $x_k.$ The solution $s_k:= \delta x_k = x_{k+1} - x_k$ is the
so-called Newton step.  Once the Newton step is obtained, then the next iterate is given by $x_{k+1}= x_k + s_k$ and 
the procedure is repeated until  convergence with respect to the prescribed tolerance is achieved. More specifically, given an 
initial approximation, say, $(v,\theta):=(v_0,v_1,\ldots, v_{N_{\xi}}, \theta_0,\theta_1,\ldots, \theta_{N_{\xi}} ) \approx (\Phi,\Lambda),$ the next approximation $(v^{+},\theta^{+})$ in the Newton's method is given by
\be
  \label{newtsys}
\left[\begin{array}{c}
v^{+}\\
\theta^{+} \\
\end{array}\right] =
\left[\begin{array}{c}
v\\
\theta \\
\end{array}\right] -
  \underbrace{\left[ \begin{array}{cc}
T(\theta) & T^{\prime}(\theta)v  \\
Q^{\prime}(v) &  0
 \end{array} \right]^{-1}}_{\mathcal{J}:=F^{\prime}}
\underbrace{\left[\begin{array}{c}
T(\theta)v\\
Q(v) \\
\end{array}\right]}_{F},
\e 
where \cite{GG07} 
\be
\label{a1}
T(\theta)= \sum_{k=0}^{N_{\xi}-1} \left[  (G_k \otimes A_k) - \theta_k (H_k \otimes  {\bf I}) \right] \in \mathbb{R}^{ N_x N_{\xi} \times N_x N_{\xi}  },
\e
\be
\label{a2}
T(\theta)v = \sum_{k=0}^{N_{\xi}-1} \left[  (G_k \otimes A_k) - \theta_k (H_k \otimes  {\bf I}) \right]v \in \mathbb{R}^{ N_x N_{\xi}   },
\e

\be
\label{a3}
T^{\prime}(\theta)v = -\sum_{k=0}^{N_{\xi}-1}  (H_k \otimes v_k) \in \mathbb{R}^{ N_x N_{\xi} \times N_{\xi}  },
\e
\be
\label{qd}{\tiny
Q(v) = {\bf d}:= \left[v^T(H_0 \otimes {\bf I})v -1, \cdots, v^T(H_{N_{\xi}-1} \otimes {\bf I})v \right]^T\in \mathbb{R}^{ N_{\xi}  },
}
\e

and
\be
\label{a4}
Q^{\prime}(v) = 2\sum_{k=0}^{N_{\xi}-1}  (H_k \otimes v_k^T) \in \mathbb{R}^{ N_x N_{\xi} \times N_{\xi}  }.
\e

\subsection{Inexact Newton method}
\label{inm}
Notwithstanding the locally quadratic convergence and simplicity of implementation of the Newton's method, 
it involves enormous computational cost, particularly when the size of the problem is  large. 
In order to  reduce the computational complexity associated with the method, Dembo, Eisenstat and Steihaug proposed 
in \cite{DES82} the {\it inexact 
Newton method} as given by Algorithm \ref{NewtonMethod}, which is a generalization of the Newton's method.

The condition in line $5$ of the algorithm is the inexact Newton condition.
Note that the real number $\eta_k$ in Algorithm \ref{NewtonMethod} is the so-called forcing term for the $k$-th 
iteration step.
At each iteration step of the inexact Newton method, $\eta_k$ should be chosen first, and then an inexact
Newton step $s_k$ is obtained by solving the Newton equations (\ref{NEq}) approximately with an efficient  solver for
 systems of linear equations. 
Quite often,  the linear system to be solved at each inexact Newton step is so large that 
it cannot be solved by direct methods. Instead,   modern iterative  solvers such as Krylov
 subspace methods  \cite{book::Saad} are typically used to solve the linear systems approximately.
This leads to a special kind of  inexact Newton method, commonly referred to as {\it inexact Newton-Krylov subspace method,}
which is  very popular in many application areas \cite{AML2007,   MK94, STW97}.  
 
We point out here that it is nevertheless hard to choose a good sequence of forcing terms.
More precisely, there may be a trade-off between the effort required to solve
the linear system to a tight tolerance and the resulting required number of nonlinear iterations. Too large a
value for $\eta_k$ results in less work for the Krylov method but more nonlinear iterations, whereas too small a
value for $\eta_k$ results in more Krylov iterations per Newton iteration. Examples of this trade-off between total
nonlinear iterations and execution time can be found in, for instance, \cite{MK94} in the context of solution of
Navier-Stokes equations. Several strategies for optimizing the computational
work with a variable forcing term $\eta_k$ are given in \cite{AML2007, EW96}. 
At any rate, it is important to note that that choice of the forcing terms should be related to specific 
problems and the information of $F(x)$ should be used
effectively \cite{AML2007}.

 \begin{algorithm}[ht]
\caption{Inexact Newton Method (INM)}
\begin{algorithmic}[1]
\small
\State Given $x_0\in \mathbb{R}^{ ({N_x}+1)N_{\xi} }$
\For{ $k=0, 1, \ldots $  ({until} $\{x_k\}$ convergence)}
\State Choose some $\eta_k\in [0,1).$
\State Solve  the Newton equations (\ref{NEq}) approximately to obtain a step $s_k$ such that
\State 
$ ||F(x_k) + F'(x_k) s_k || \leq \eta_k ||F(x_k)||.$
\State Set  $x_{k+1}= x_{k} + s_k. $
\EndFor
\end{algorithmic}
\label{NewtonMethod}
\end{algorithm}
For practical computations, there are
 some concrete strategies, one of  which was proposed originally by Dembo and Steihaug in \cite{DS83}, namely, 
\be
\label{eta1}
\eta_k = \min\{1/(k+2), ||F(x_k)|| \}.
\e
Moreover, Cai et. al in  \cite{Cai95} propose the following constant forcing terms:
\be
\label{eta2}
\eta_k = 10^{-4}.
\e
Two other popular adaptive strategies were proposed by Eisenstat and Walker in \cite{EW96}: 
\begin{enumerate}
 \item[(a)] Given some $\eta_0\in [0,1),$ choose
\be
 \label{eta3}
 \eta_k = \begin{cases}
   \zeta_k, \quad\quad\quad \quad\quad\quad\quad\quad\quad    \eta^{(1+\sqrt{5})/2}_{k-1}\leq 0.1,\nonumber\\
  \max\left\{\zeta_k, \eta^{(1+\sqrt{5})/2}_{k-1}\right\},  \quad \eta^{(1+\sqrt{5})/2}_{k-1}> 0.1,
 \end{cases}
\e 
where
\[
\zeta_k = \frac{||F(x_k) - F(x_{k-1})  - F'(x_{k-1}) s_{k-1} ||}{||F(x_{k-1})||}, \quad k=1,2\ldots
\]
or
\[
\zeta_k = \frac{|\left|\left|F(x_k)\right|\right| - \left|\left|F(x_{k-1})  + F'(x_{k-1}) s_{k-1} \right|\right||}{||F(x_{k-1})||}, \quad k=1,2\ldots
\]
 \item[(b)] Given some $\tau\in [0,1), \; \omega\in [1,2),\; \eta_0\in [0,1),$ choose
\be
 \label{eta4}
 \eta_k = \begin{cases}
   \zeta_k, \quad\quad\quad \quad\quad\quad\quad\quad\quad    \tau \eta^{\omega}_{k-1}\leq 0.1,\nonumber\\
  \max\left\{\zeta_k, \tau \eta^{\omega}_{k-1}\right\},  \quad\quad\;\;\; \tau \eta^{\omega}_{k-1}> 0.1,
 \end{cases}
\e 
where
\[
\zeta_k = \tau\left(\frac{||F(x_k)||}{||F(x_{k-1})||}\right)^{\omega}, \quad k=1,2\ldots
\]
\end{enumerate}

The numerical experiments in \cite{EW96} show that the above two choices $(a)$ and $(b)$ can 
effectively overcome the `over-solving'
phenomenon, and thus improve the efficiency of the inexact 
Newton method\footnote{The concept of `over-solving' implies that at early Newton iterations $\eta_k$ is too small. Then one may obtain
an accurate linear solution to an inaccurate Newton correction. This may result in a poor Newton update
and degradation in the Newton convergence. In \cite{TWS02} it has been demonstrated that in some situations
the Newton convergence may actually suffer if $\eta_k$ is too small in early Newton iterations.
}. 
In particular,
the authors added  safeguards (bounds) to choice (a) and (b)
to prevent the forcing terms from becoming too small too quickly, so that
 more concrete strategies are obtained.
Besides, choice (a) and choice (b) with
$\tau\geq 0.9$ and $\omega\geq (1 +\sqrt{5})/2$ have the best performances. We adopt choice (b) in our numerical experiments.

 The inexact Newton method is locally convergent as shown in the following result from \cite{DES82}.
\begin{theorem}\cite[Theorem 2.3]{DES82}
\label{thminm}
Assume that $F : \mathbb{R}^n \rightarrow \mathbb{R}^n$
 is continuously differentiable, $x^{\ast}\in \mathbb{R}^n$
such that $F(x^{\ast}) = 0$ and $F^{\prime}(x^{\ast})$
is nonsingular. Let $0<\eta_k<\eta_{\max} <t <1$ be  given constants. 
If the forcing terms $\{\eta_k\}$
in the inexact Newton method satisfy $\eta_k \leq \eta_{\max} < t <1$ for all $k$, then there exists $\varepsilon>0,$ such that for any 
$x_0\in N_{\varepsilon}(x^{\ast}):= \{x : ||x - x^{\ast}|| < \varepsilon \}$, the sequence $\{x_k\}$ 
generated by the inexact Newton method converges to $x^{\ast}$
, and
\[
|| x_{k+1} - x_k ||_{\ast}\leq t || x -  x^{\ast} ||_{\ast},
\]
where $||y||_{\ast} =||F^{\prime}(x^{\ast})y||.$
\end{theorem}

By Theorem \ref{thminm}, if the forcing terms $\{\eta_k\}$ in the inexact Newton method are uniformly strict less than $1$,
then the method is locally convergent. The convergence rate of the inexact Newton method is, moreover,  established in
the following result from \cite{DES82}.

\begin{theorem}\cite[Corollary 3.5]{DES82}
\label{thminm2}
Assume that $F : \mathbb{R}^n \rightarrow \mathbb{R}^n$
 is continuously differentiable, $x^{\ast}\in \mathbb{R}^n$
such that $F(x^{\ast}) = 0$ and $F^{\prime}(x^{\ast})$
is nonsingular. 
If the sequence $\{x_k\}$ generated by inexact Newton method converges
to $ x^{\ast},$ then
\begin{itemize}
 \item $\{x_k\} \rightarrow  x^{\ast}$ super-linearly when $\eta_k  \rightarrow 0.$
 \item  $\{x_k\} \rightarrow  x^{\ast}$ quadratically when $\eta_k = \mathcal{O}(||F(x_k)||)$
 and $||F^{\prime}(x) ||$ is Lipschitz continuous at $x^{\ast}$.
\end{itemize}
\end{theorem}

For more details of local convergence theory and the
role played by the forcing terms in inexact Newton methods, see e.g., \cite{AML2007, EW96}. 
We proceed next to give an overview of Krylov subspace methods.


\subsection{Krylov subspace methods}
Krylov subspace methods  are probably  the most popular methods for solving large, 
sparse linear systems (see e.g. \cite{ESW14} and the references therein).
The basic idea behind  Krylov subspace methods is the following.
 Consider, for arbitrary $A\in \mathbb{R}^{m \times m}$ and $b\in \mathbb{R}^{m},$  the linear system
\be
\label{testsys}
Ax = b.
\e
Suppose now that $x_0$ is an initial guess for the solution $x$ of (\ref{testsys}), and define
the initial residual $r_0 = b - Ax_0.$ Krylov subspace methods are iterative
methods whose $k$th iterate $x_k$ satisfies
\footnote{Krylov  methods require only matrix-vector products to carry out the iteration (not the
individual elements of A) and this is key to their use with the Newton's method, as will be seen below.}
\be
\label{krylov}
x_k  \in x_0 + \mathbb{K}_k(A,x_0), \;\;  k=1,2,\ldots,
\e
where
\be
\mathbb{K}_k(A,x_0):= \mbox{span}
\left\{r_0,Ar_0 ,\ldots,A^{k-1}r_0 \right\}
\e
denotes the $k$th Krylov subspace generated by $A$ and $r_0.$ The
 Krylov subspaces form a nested sequence that ends with dimension
$d = \mbox{dim}( \mathbb{K}_{m}(A,r_0 ) )\leq m,$ i.e.,
\[
\mathbb{K}_{1}(A,r_0 )\subset \ldots \subset \mathbb{K}_{d}(A,r_0 ) =\cdots =  \mathbb{K}_{m}(A,r_0 ).
\]
In particular, for each $k \leq d,$ the Krylov subspace $\mathbb{K}_k(A,r_0 )$ has dimension
$k.$ Because of the $k$ degrees of freedom in the choice of the iterate $x_k,$ $k$
constraints are required to make $x_k$ unique. In Krylov subspace methods
this is achieved by requiring that the $k$th residual $r_k =  b - Ax_k$ is orthogonal (with respect to the Euclidean inner product)
to a $k$-dimensional space $\mathcal{C}_k,$ called the constraints space:
\be
\label{cspace}
r_k = b - Ax_k \in r_0 + A\mathbb{K}_k(A,r_0 ),
 \e
where $r_k \perp \mathcal{C}_k.$
It can be shown \cite{BGL05} that there exists a uniquely defined iterate $x_k$ of the form (\ref{krylov}) and for which
the residual $r_k = b - Ax_k$ satisfies (\ref{cspace}) if
\begin{itemize}
 \item [(a)] $A$ is symmetric positive definite  and $\mathcal{C}_k=\mathbb{K}_k(A,r_0 ),$ or
\item  [(b)] $A$  is nonsingular and $\mathcal{C}_k=A\mathbb{K}_k(A,r_0 ).$
\end{itemize}
In particular, (a)  characterizes the conjugate gradient (CG) method \cite{ESW14}  whereas (b) characterizes the minimal
residual (MINRES) method \cite{minres}, the generalized minimal residual
(GMRES) method \cite{gmres},  and the  bi-conjugate gradient stabilized (BiCGstab)  method \cite{Vorst92}.

A vast majority of fully coupled 
nonlinear applications of primary interest (including the one considered herein)
result in Jacobian matrices that are non-symmetric. A
further point of discrimination is whether the method is derived from the long-recurrence Arnoldi
orthogonalization procedure, which generates orthonormal bases of the Krylov subspace, or the short-recurrence
Lanczos bi-orthogonalization procedure, which generates non-orthogonal bases for non-symmetric matrices $A$.

Note that GMRES is an Arnoldi-based method.
In GMRES, the Arnoldi basis vectors form the trial subspace out of which the solution is constructed. One
matrix-vector product is required per iteration to create each new trial vector, and the iterations are terminated
based on a by-product estimate of the residual that does not require explicit construction of
intermediate residual vectors or solutions -- a major beneficial feature of the algorithm. GMRES has a
residual minimization property in the Euclidean norm (easily adaptable to any inner-product norm) but
requires the storage of all previous Arnoldi basis vectors. Full restarts, seeded restarts, and moving fixed sized
windows of Arnoldi basis vectors are all options for fixed-storage versions. Full restart is simple and
historically the most popular, though seeded restarts show promise. The 
BiCGstab  methods \cite{Vorst92} are Lanczos-based alternatives to GMRES for non-symmetric problems.
In BiCGstab methods, the Lanczos basis
vectors are normalized, and two matrix-vector products are  required per iteration. However, these methods
enjoy a short recurrence relation, so there is no requirement to store many Lanczos basis vectors. These
methods do not guarantee monotonically decreasing residuals.
We refer to \cite{ESW14} for more details on Krylov methods, and for preconditioning for linear problems.

\subsection{Inexact Newton-Krylov method with backtracking}
In practice, {\it globalization
strategies} leading from a convenient initial iterate into the ball of convergence of Newton's  method around
the desired root are often required to enhance the robustness of the inexact Newton method.
More precisely, globalization implies 
 augmenting  Newton's  method with certain auxiliary procedures  that increase
the likelihood of convergence when good initial approximate solutions are not available.
Newton-Krylov methods, like all Newton-like methods, must usually be  globalized.
Globalizations are typically structured to test whether a step gives satisfactory
progress towards a solution and, if necessary, to modify it to obtain a step that does 
give satisfactory progress \cite{PSSW06}. A major class of  globalization approaches\footnote{See e.g. \cite{PSSW06, STW97} 
for a detailed discussion on other  globalization strategies such as trust-region methods.} which we consider in this 
paper are the {\it backtracking
(line-search, damping) methods}. In these methods, the step lengths are adjusted (usually
shortened) to obtain satisfactory steps.
On the one hand,  backtracking methods have  the attrative feature of  the relative ease
with which they can be implemented; on the other hand,
each step direction in these methods is restricted to be that of the initial trial step, which may
be a weak descent direction, especially if the Jacobian is ill-conditioned \cite{STW97}.

The inexact Newton backtracking method (INBM) is given in Algorithm~\ref{INB}. In this algorithm,
the backtracking globalization resides in the while-loop, in
which steps are tested and shortened as necessary until the acceptability condition
\be
\label{ac}
||F(x_k + s_k )|| \leq [1 - t(1 - \eta_k)] ||F(x_k)||, 
\e
holds.  As noted in \cite{EW94}, if $F$ is continuously differentiable, then this globalization
produces a step for which (\ref{ac}) holds after a finite number of passes through the
while-loop; furthermore, the inexact Newton condition (cf. line $5$ in Algorithm \ref{NewtonMethod}) still holds for the final
$s_k$ and $\eta_k$. The condition (\ref{ac}) is a `sufficient-decrease' condition on
$||F(x_k + s_k )||.$

\begin{algorithm}[t]
\caption{Inexact Newton Backtracking Method (INBM)}
\label{INB}
\begin{algorithmic}[1]
\small
\State Let  $x_0\in \mathbb{R}^{ ({N_x}+1)N_{\xi}}, \; \eta_{\max}\in [0,1), \; t\in (0,1),\; \mbox{and}\; 0<\theta_{\min} <\theta_{\max} <1,\;$ be given.
\For{ $k=0, 1, \ldots $ ({until} $\{x_k\}$ convergence) }
\State  Choose initial ${\eta}_k\in [0,\eta_{\max})$ and
solve (\ref{NEq}) approximately to obtain $s_k$ such that
\State $\label{necon2} ||F(x_k) + F'(x_k) s_k || \leq {\eta}_k ||F(x_k)||.$
\State While {$ ||F(x_k + s_k )|| > [1 - t(1 - \eta_k)] ||F(x_k)|| $ }
\State   Choose  $\theta\in [\theta_{\min}, \theta_{\max}]. $
\State   Update $s_k \leftarrow \theta s_k$  and $\eta_k \leftarrow 1 -\theta(1-\eta_k).$
\State  Set  $x_{k+1}= x_{k} + s_k. $
\EndFor
\end{algorithmic}
\label{Newtonback}
\end{algorithm}

In \cite{EW96}, the  authors  
show with experiments that backtracking globalization  significantly improves the robustness of a Newton-GMRES
method when applied to nonlinear problems, especially when combined with
adaptively determined forcing terms. In this work, we combine 
the backtracking globalization with low-rank techniques to tackle the high-dimensional 
stochastic eigenproblem. Our motivation for employing low-rank techniques stems from the fact that
despite the  advantages of the INKM with backtracking in solving nonlinear problems,
 for the stochastic problem (\ref{randeqn}) -- (\ref{normeig}) under consideration, the dimensions of the  Jacobian  quickly become prohibitively
 large with respect to the discretization parameters.
 As a consequence, one expects overwhelming memory and computational time requirements,
 as the block-sizes of the Jacobian matrix become vast. 
 This is a major drawback of the SGM.
 In this paper, we propose to tackle this {\it curse of dimensionality} 
 with a low-rank version of INKM. Low-rank strategies have proven to be quite efficient 
 in solving problems of really high computational complexity arising, for instance, from 
 deterministic and stochastic time-dependent optimal control problems \cite{BSOS215, BOS15a,  StoB13}, 
 PDEs with random coefficients \cite{BOS14, AO16}, etc. 
 The low-rank technique presented here only needs to store a small portion of the vectors in comparison
 to the full problem and  we want present
  this approach in the sequel.
  
\subsection{Low-rank inexact Newton-Krylov method}
As we have already noted earlier, we will use a Krylov solver algorithm  as an optimal solver for
the Newton equation (cf. (\ref{NEq}) in step $3$  in Algorithm~\ref{INB})
in each INKM iteration. 
In particular, our approach is based on the low-rank version of the chosen Krylov solver. 
Although the low-rank method discussed herein can be easily extended to other Krylov solvers
\cite{BOS14, BOS15a,  StoB13}, we focus mainly on
 BiCGstab \cite{KT11}. In this section, 
we proceed first to  give a brief overview of this  low-rank iterative solver. Now, recall first that 
\be
 \label{kronvec}
 \mathrm{vec}(WXV) = (V^T\otimes W)\mathrm{vec}(X),
\e
where    $\mathrm{vec}(X) = (x_1,\ldots,x_p)^T \in \mathbb{R}^{np \times 1}$  is a column vector obtained by stacking the columns of 
 the matrix $X=[x_1,\ldots,x_p]\in \mathbb{R}^{n\times p}$ on top of each other. Observe then that, using (\ref{kronvec}),  each  Newton equation
 (\ref{NEq}) can be rewritten  as $ \mathcal{J}  \mathcal{X}= \mathcal{R},$ where
\[
 \mathcal{J}:=F^{\prime}  =   \left[\begin{array}{cc}
\sum\limits_{i=0}^{N_{\xi}-1} \left[(G_i - \lambda_iH_i) \otimes (A_i - I_{N_x})\right]  &  -\sum\limits_{i=0}^{N_{\xi}-1} H_i \otimes v_i \\
2\sum\limits_{i=0}^{N_{\xi}-1} H_i \otimes v^T_i & 0 \\
\end{array}\right],
\]
\[
\mathcal{X}:=s= \left[\begin{array}{c}
{\mbox{vec} (Y)} \\
{\mbox{vec} (Z)} \\
\end{array}\right], 
\;\;\;
\mathcal{R}:=-F=\left[\begin{array}{c}
{\mbox{vec}(R_1)} \\
{\mbox{vec}(R_2)} \\
\end{array}\right],
\]
and 
\[
R_1=\mbox{vec}^{-1}\left(\sum\limits_{i=0}^{N_{\xi}-1} \left[(G_i - \lambda_iH_i) \otimes (A_i - I_{N_x})\right]v\right),\;\; R_2={\mbox{vec}^{-1}(\bf d)},
\]
where ${\bf d}$ is as given by (\ref{qd}).
Hence, (\ref{kronvec})  implies that
{\small{
\be
\label{rhseq}
 \mathcal{J}  \mathcal{X}= \mbox{vec}\left( \sum\limits_{i=0}^{N_{\xi}-1} \left[\begin{array}{c}
 (A_i - I_{N_x}) Y (G_i - \lambda_iH_i)^T  - v_i Z H_i^T  \\
  2v^T_i Y H_i^T\\
\end{array}\right]\right) = \mbox{vec}\left( \left[\begin{array}{c}
{R_1} \\
{R_2} \\
\end{array}\right]\right).
\e
}}

Our approach is essentially based on the assumption  
that both the solution matrix $\mathcal{X}$  admits a low-rank representation; that is,
\be
 \label{lowrankapp}
 \begin{cases}
   Y = W_YV_Y^T, \;\;\mbox{with} \;\;W_Y\in \mathbb{R}^{(N_x +1)\times k_1},\;\; V_Y\in \mathbb{R}^{N_{\xi}\times k_1}\\
   Z = W_ZV_Z^T, \;\;\mbox{with} \;\;W_Z\in \mathbb{R}^{(N_x +1)\times k_2},\;\; V_Z\in \mathbb{R}^{N_{\xi}\times k_2}
\end{cases}
\e
where $k_{1,2,3}$ are small relative to $N_{\xi}.$ Substituting (\ref{lowrankapp}) in (\ref{rhseq}) and ignoring the vec  operator,
we then obtain\footnote{Note that $v_i$ in (\ref{newrhs})  comes from the previous low-rank iterate of the nonlinear Newton solver.}
\be
\label{newrhs}
 \sum\limits_{i=0}^{N_{\xi}-1} \left[\begin{array}{c}
 (A_i - I_{N_x}) W_YV^T_Y (G_i - \lambda_iH_i)^T  - v_i W_ZV^T_Z H_i^T\\
  2v^T_i W_YV^T_Y H_i^T\\
\end{array}\right]=
\left[\begin{array}{c}
{R_{11}R_{12}^T } \\
{R_{21}R_{22}^T} \\
\end{array}\right],
\e
where $R_{11}R_{12}^T $  and $R_{21}R_{22}^T $  are the low-rank representations of  
$R_1$  and $R_2,$  respectively.

The attractiveness of this approach lies therefore in the fact that one can rewrite the three block rows
in the left hand side  in (\ref{newrhs}), respectively, as 
\be
{\small{
 \label{lrformat}
 \begin{cases}
\mbox{ (first block row)}  \sum\limits_{i=0}^{N_{\xi}-1} \left[\;\;  (A_i - I)W_Y  \;\;\;\;- v_iW_Z \;\;\right]  \left[ \begin{array}{ll}
          V_Y^T(G_i - \lambda_iH_i)^T\\
       V_Z^TH_i^T      \end{array} \right],\\\\
  \mbox{ (second block row)} \sum\limits_{i=0}^{N_{\xi}-1} \left[\;\;  2v_iW_Y   \;\;\right]  \left[ \begin{array}{ll}
          V_Y^TH_i^T\\
           \end{array} \right],
\end{cases}
}}
\e
so that the low-rank nature of the factors guarantees fewer multiplications with the submatrices  while maintaining smaller storage requirements.
More precisely, keeping in mind that
\begin{equation}
\label{vertx}
x= \mbox{vec}\left( \left[\begin{array}{c}
{X_{11}X_{12}^T } \\
{X_{21}X_{22}^T } \\
\end{array}\right]\right)
\end{equation}
corresponds to the associated vector $x$ from a vector-based version of the  Krylov solver,  matrix-vector multiplication 
in our low-rank Krylov solver is given by Algorithm~\ref{matvec}.

{\small{
\begin{algorithm}[t]
  \caption{Jacobian-vector multiplication in low-rank format $\tt{Amult}$}
\label{matvec}
  \begin{algorithmic}[1]
\State Input: $W_{11},W_{12},W_{21},W_{22}$
\State Output: $X_{11},X_{12},X_{21},X_{22}$
\State $X_{11}=\sum\limits_{i=0}^{N_{\xi}-1} \left[\;\;  (A_i - I)W_{11}  \;\;\;\;- v_iW_{21} \;\; \right]$
\State $X_{12}= \left[\;\;  (G_i - \lambda_{i}H_i) W_{12}  \;\; \cdots \;\;  H_iW_{22} \;\;\right],\;\; i=0,\cdots,N_{\xi}-1.$
\State $X_{21}=\sum\limits_{i=0}^{N_{\xi}-1}\left[\;\;   2v^T_iW_{11}   \;\; \right]$
\State $X_{22}= \left[ \;\; H_i W_{12}  \;\; \right], \;\; i=0,\cdots,N_{\xi}-1.$
 \end{algorithmic}
\end{algorithm}
}}

Note that an important feature of  low-rank Krylov solvers is that  the iterates of the  solution
matrices $Y$ and $Z$ in the algorithm are truncated 
by a truncation operator $\mathcal{T}_\epsilon$ with a prescribed tolerance $\epsilon.$  This could be accomplished via QR decomposition 
as in \cite{KT11}  or truncated singular value decomposition  (SVD)  as in \cite{BOS14, StoB13}.
The truncation operation is necessary because the new computed factors could have  increased ranks compared to the original
factors in (\ref{lrformat}). Hence, a truncation of all the factors after the matrix-vector
products, is used to construct new
factors; for instance,
{\small{\[
[\tilde{X}_{11}, \tilde{X}_{12}]:= \mathcal{T}_\epsilon\left( [{X}_{11}, {X}_{12}] \right)
=\mathcal{T}_\epsilon\left(\sum\limits_{i=0}^{N_{\xi}-1} \left[\;\;  (A_i - I)W_Y  \;\;\;\;- v_iW_Z \;\;\right]  
\left[ \begin{array}{ll}
          V_Y^T(G_i - \lambda_iH_i)^T\\
       V_Z^TH_i^T      \end{array} \right]
       \right).
\] } }
Moreover, in order to  ensure that the inner products within the iterative low-rank solver are  computed efficiently, we use the fact that 
\[
\left<x,y\right>=\vec{X}^T\vec{Y}=\trace{X^TY}
\]
to deduce that
\begin{eqnarray}
\label{tracexy}
\trace{X^TY} &=& \trace{\underbrace{\left(X_{11}X_{12}^{T}\right)^{T}}_{\textrm{Large}}\underbrace{\left(Y_{11}Y_{12}^{T}\right)}_{\textrm{Large}}
 + \underbrace{\left(X_{21}X_{22}^{T}\right)^{T}}_{\textrm{Large}}\underbrace{\left(Y_{21}Y_{22}^{T}\right)}_{\textrm{Large}}  } \nonumber\\
 &=& \trace{\underbrace{Y_{12}^{T}X_{12}}_{\textrm{Small}}\underbrace{X_{11}^{T}Y_{11}}_{\textrm{Small}}
+ \underbrace{Y_{22}^{T}X_{22}}_{\textrm{Small}}\underbrace{X_{21}^{T}Y_{11}}_{\textrm{Small}}  },
\end{eqnarray}
where $X$ and $Y$ are as given in (\ref{vertx}), which allows us to  compute the trace of small matrices rather than of the ones
from the full model.

For more details on implementation issues, we refer the interested reader to \cite{BOS14, StoB13}.

\subsection{Preconditioning}
The purpose of preconditioning the INBM is to reduce the number of Krylov iterations,
as manifested  by efficiently clustering eigenvalues of
the iteration matrix. Traditionally, for linear problems, one chooses a few iterations of a simple iterative
method (applied to the system matrix) as a preconditioner. 
Throughout this paper, we will focus  mainly on  mean-based block-diagonal preconditioners.
More specifically,  we precondition the Jacobian matrix $J$ (cf.  (\ref{newtsys}) )
in  the INBM algorithm   with a preconditioner $\mathcal{P}$ of the form
\be 
\mathcal{P} :=  
  \left[ \begin{array}{cc}
E & 0  \\
0 & S 
 \end{array} \right],\label{idealp}
\e 
where
\be
\label{schurcomp}
S=CE^{-1}B
\e
is the (negative) {\it Schur complement}. Moreover, $E:= T(\Lambda), \; B:= T^{\prime}(\Lambda)$ and 
$C:=Q^{\prime}(\Phi)$ as given, respectively, by  (\ref{a1}), (\ref{a3}) and (\ref{a4}).
We note here that (\ref{idealp}) is only an ideal preconditioner 
for the Jacobian
in the sense that it is not cheap to solve the system with it. 
In practice, one often has to approximate its two diagonal blocks in order to use $\mathcal{P} $ with Krylov solvers. 
Here, we propose to  approximate the $(1,1)$ blocks with $(G_0 - \lambda_0H_0) \otimes (A_0 - I_{N_x})$ which is easy to invert:
if we use the normalized Legendre polynomial chaos to compute the matrices $G_i$ and $H_i$, then 
$(G_0 - \lambda_0H_0)= (1 - \lambda_0)I_{N_{\xi}}$ so that action of the approximated $(1,1)$ block is just $N_{\xi}$ copies of $(A_0 - I_{N_x}).$
 To approximate the Schur complement $S,$ 
 that is, block $(2,2),$ poses more difficulty, however. One possibility 
is to approximate  $S$ by dropping all but the first terms in $B, C$  and $E$ to obtain
\begin{eqnarray}
{S}_0 &:=& 2(1 - \lambda_0)^{-1}(I_{N_{\xi}} \otimes v_0)( I_{N_{\xi}}\otimes (A_0 - I_{N_x})^{-1}) (I_{N_{\xi}}\otimes v_0)^T \nonumber\\
 &=& 2(1 - \lambda_0)^{-1}I_{N_{\xi}}\otimes \left[ v_0(A_0 - I_{N_x})^{-1}v^T_0  \right].\label{shur1}
\end{eqnarray}
This is the version we use in our experiments, and its
implementation  details are provided in Algorithm \ref{aprec}.

{\small{
\begin{algorithm}[t]
  \caption{Preconditioner implementation in low-rank Krylov solver }
\label{aprec}
  \begin{algorithmic}[1]
\State Input:  $W_{11},W_{12},W_{21},W_{22}$
\State Output: $X_{11},X_{12},X_{21},X_{22}$
\State Solve: $(A_0 - I_{N_x}) X_{11}=W_{11} $ 
\State Solve: $(1 - \lambda_0)X_{12}=W_{12}$
\State Solve: $\left[ v_0(A_0 - I_{N_x})^{-1}v^T_0  \right] X_{21}=W_{21} $ 
\State Solve: $2(1 - \lambda_0)^{-1}X_{22}=W_{12}$
 \end{algorithmic}
\end{algorithm}
}}

\section{Numerical results}
\label{numresults}
In this section, we present some numerical results
obtained with the proposed inexact Newton-Krylov solver for the stochastic eigenproblems (\ref{randeqn}).
 The numerical experiments were performed on a Linux machine with 80 GB RAM using MATLAB\textsuperscript{\textregistered}
7.14 together with a MATLAB version of the 
algebraic multigrid (AMG) code HSL MI20 \cite{BMS07}.
 We implement  our mean-based preconditioner  using one 
V-cycle of AMG with symmetric Gauss-Seidel (SGS) smoothing to approximately invert $A_0 - I_{N_x}.$ 
We remark here that we apply the method as a black-box in each experiment and the set-up of the approximation to $A_0 - I_{N_x}$ only needs to be performed once.
Unless otherwise stated, in all the simulations, BiCGstab is terminated when the relative residual error is reduced to $tol = 10^{-5}.$ Note that  $tol$ 
should be chosen such that the truncation tolerance $trunctol \leq tol;$ 
otherwise, one would be essentially iterating on the `noise' from the low-rank truncations. In particular, we have chosen herein  $trunctol = 10^{-6}.$
We have used the Frobenius  norm throughout our numerical experiments.

Before proceeding to present our numerical example, it is perhaps pertinent to highlight certain factors that often influence the convergence of
 the inexact Newton  method \cite{FBF15}:
\begin{itemize}
 \item  the proximity of the initial guess. Here, we have employed uniformly distributed samples for our initial guess.
 \item  The globalization technique employed, (e.g. backtracking,
or trust region). In this paper, we have used only backtracking  and it worked quite well for our considered problem.
 \item The discretization of the SEVPs -- failure of the spatial discretization to
adequately reflect the underlying physics of the continuous problem can cause
convergence difficulties for globalized Newton-Krylov methods.
\item The convergence of the Krylov solver and preconditioning strategy employed -- 
using nonlinear preconditioning techniques can be an alternative \cite{BKST13}.
\end{itemize}

For our numerical example, let  $\mathcal{D}=(0,1)\times  (0,1)$. We consider  the stochastic eigenproblem of finding the functions
$\lambda: \Omega \rightarrow \mathbb{R}$ and $\varphi: \Omega \times D \rightarrow  \mathbb{R}$ 
such that, $\mathbb{P}$-almost surely, the following   holds:
\be
 \label{modelprob}
\left\{
\begin{aligned}
     - \nabla\cdot(a({\cdot},\omega)\nabla \varphi({\cdot},\omega)) &= \lambda(\omega)\varphi({\cdot},\omega),
       \;\;\mbox{in} \; \;\mathcal{D}\times  \Omega,  \quad\\
    \varphi({\cdot},\omega) &= 0, \;\;\mbox{on} \;\; \partial\mathcal{D}\times  \Omega, \quad\\
    \end{aligned}
\right.
\e
where $a:\mathcal{D}\times \Omega \rightarrow \mathbb{R}$ is a random coefficient field. 
 We assume that 
there   exist positive constants $a_{\min}$ and $a_{\max}$ such that
\be
\label{pal}
\mathbb{P}\left( \omega\in \Omega:a({\bf x},\omega)\in [a_{\min}, a_{\max}],\, \forall {\bf x}\in\mathcal{D}  \right)=1.
\e
Here,  the random input   $a(\cdot,\omega)$ admits a KLE and has a covariance function given by  
 \[
 C_a({\bf x},{\bf y})= \sigma_a^2\exp\left(-\frac{|x_1-y_1|}{\ell_1}  -\frac{|x_2-y_2|}{\ell_2}\right), \;\;\;\forall ({\bf x},{\bf y}) \in [-1,1]^2,
\]
with correlation lengths $\ell_1=\ell_2=1 $  and  mean of the random field  $a$  in the model $\mathbb{E}[a]=1$.
 The forward problem has been extensively studied  in, for instance, \cite{PE09}. 
The eigenpairs of the KLE of the random field $a$ are given explicitly in \cite{GS91}.
Note then that discretising in space yields the expression (\ref{randeqn}) with the random matrix
$\mathcal{A}(\omega)$ having the explicit expression
(\ref{rmatrix}). 
In particular, the stiffness matrices $A_k\in \mathbb{R}^{N_x\times N_x}, \; k=0,1,\ldots, m,$ in (\ref{rmatrix}) are given, respectively, by
 \begin{alignat}{2}
A_0(i,j)  &= \int_{\mathcal{D}} \mathbb{E}[a](x)\nabla \phi_i(x)\nabla \phi_j(x)\; dx, & \label{K0}\\
A_k(i,j)  &= \sigma_a \sqrt{\gamma_k} \int_{\mathcal{D}} \vartheta_k(x)\nabla \phi_i(x)\nabla \phi_j(x)\; dx, &\label{Kj} \;\; k>0,
\end{alignat}
where $\sigma_a$ is the standard deviation  of $a$.
Here, $\{\gamma_k\}$ and  $\{ \vartheta_k(x)\}$ are, respectively, the  eigenvalues and
eigenfunctions corresponding to a covariance function associated with $a.$
Also, $\{\phi_j(x)\}$  are  ${\bf Q}_1$ spectral elements which we have used to discretize the problem (\ref{modelprob}) in the spatial domain $\mathcal{D}.$
Moreover, we choose $\xi=\{\xi_1,\ldots,\xi_m\}$ such that $\xi_k\sim \mathcal{U}[-1,1],$ and $\{\psi_k\}$ are $m$-dimensional 
Legendre polynomials with support in $[-1,1]^m.$  
In particular, we have $N_{\xi}=210$ (with $m=6$ and $r=4$; cf. (\ref{Nxi}) ).

In what follows, we consider two cases.
First, in Table \ref{bicgstab_eta}, we set $\sigma_a=0.01$ and  $N_x =49,$ so that from (\ref{Jsize}) and  (\ref{newtsys}), we have a Jacobian matrix $\mathcal{J}$  of dimension
dim$(\mathcal{J}):=({N_x}+1)N_{\xi}=10,500.$  
Admittedly,
one would argue that this dimension of the Jacobian is small, and as such can as well be handled without the low-rank method proposed in this paper!
Such an arguement is understandably valid. However, this is basically intended to provide a first and simple insight as to how the algorithm works. A more difficult case is provided
in Table \ref{bicgstab_eta2} where we have increased  $\sigma_a$ and $N_x$ to  $\sigma_a=0.1$ and $N_x=392,704,$ respectively. Hence, we obtain a Jacobian of size
dim$(\mathcal{J}):=({N_x}+1)N_{\xi}=82,468,050$ at each 
inexact Newton iteration! We note here that,   besides the increased dimension, increasing
the variability ($\sigma_a$ ) in the problem equally
increases the complexity of the linear system to be solved \cite{BOS14}. 

The methods discussed in the previous sections have many parameters that must
be set, e.g., the maximum BiCGstab iterations, maximum forcing term, etc. These parameters
affect the performance of the methods. We chose parameters commonly used in the
literature.
In particular, for the forcing terms $\eta_{k}$, we
 set $\eta_{0} = 0.9, \eta_{\max} = 0.9, \; \eta_{\min} = 0.1.$
For the  backtracking parameters, we  used $\theta_{\max} = 0.1, \; \theta_{\max} = 0.5$; 
the maximum number of backtracks allowed is $20.$

Now, we consider the first case; that is, when dim$(\mathcal{J}):=({N_x}+1)N_{\xi}=10,500.$ We note here that the INBM algorithm  presented in this paper computes one eigenvalue nearest to the initial
guess. To compute two or more distinct (or multiple) roots of $F(x)=0$ for an SEVP would require some specialized techniques, 
which  can be a lot more involved. Nevertheless, this task is currently under investigation, and
 an algorithm for the computation of other eigenvalues will be presented in a subsequent 
paper. 

\begin{figure}[bt]
\centering
 \caption{Eigenvalues of the deterministic matrix $A_0$}
\includegraphics[width=0.85\textwidth,height=0.5\textwidth]{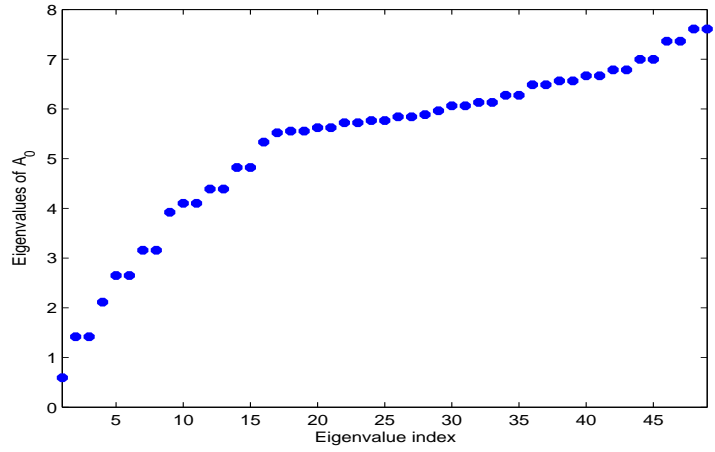}
\label{eigs_A0}
\end{figure}

All the eigenvalues of the deterministic matrix (i.e. $A_0$) are shown in Figure~\ref{eigs_A0}.
The first six smallest eigenvalues of $A_0$ are $ 0.5935, 1.4166, 1.4166, 2.1143, 2.6484, 2.6484.$
 We note here that most of the eigenvalues of the matrix $A_0$  are  either repeated
or quite clustered.
Observe in particular that $1.4166$  and $2.6484$ are repeated eigenvalues.

\begin{figure}[t]
\centering
 \caption{Convergence of low-rank INBM for the second stochastic eigenvalue $\lambda_2(\omega).$}
\includegraphics[width=0.85\textwidth,height=0.51\textwidth]{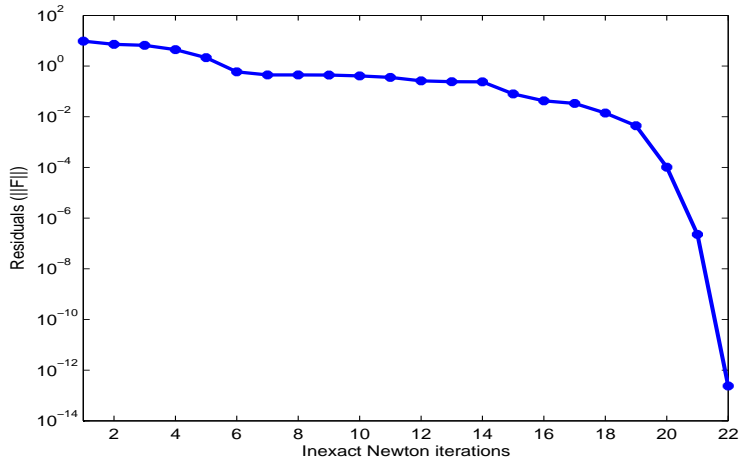}
\label{conv_plot}
\end{figure}

\begin{figure}[h]
\centering
 \caption{Probability density function (pdf) estimate of the second eigenvalue obtained with $\sigma_a=0.1$ }
\includegraphics[width=0.85\textwidth,height=0.51\textwidth]{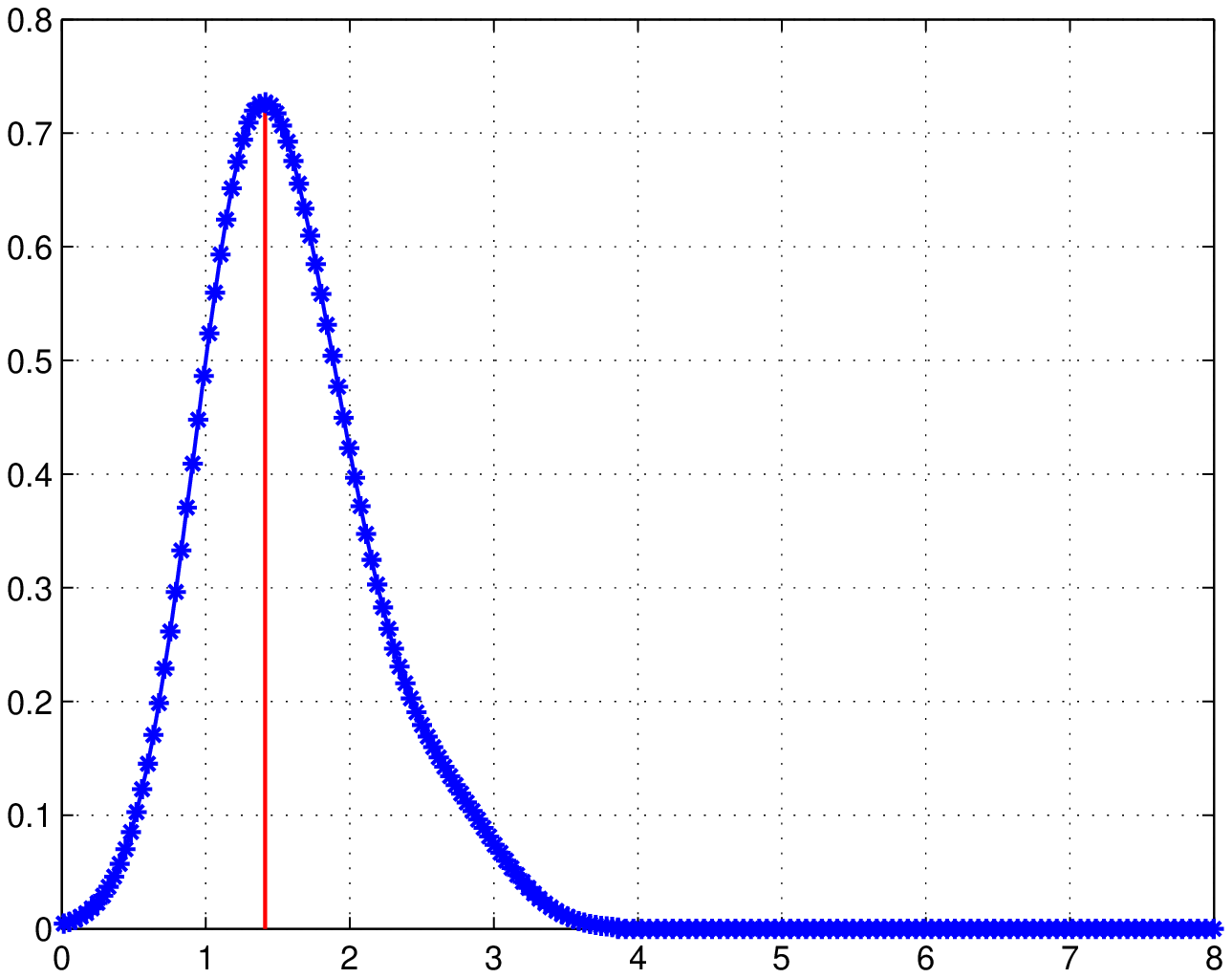}
\label{pdf_plot_2eig}
\end{figure}

In Figure~\ref{conv_plot}, we show the convergence of the low-rank INBM 
to  the second stochastic eigenvalue  $\lambda_2(\omega).$ The figure confirms the super-linear convergence
of the inexact Newton algorithm as we reported earlier.
In  Table \ref{tractab} and Figure \ref{pdf_plot_2eig}, we have shown  the
first eight coefficients  of the spectral expansion gPCE   and the probability density function (pdf) of the second eigenvalue, respectively.
Observe here that the pdf is as expected centered at the mean of the stochastic eigenvalue, i.e $1.4121.$
 This quantity can also be obtained from the first coefficient of the gPCE in Table \ref{tractab}, since
 from (\ref{lambdas}), we have
 \[
\mathbb{E}(\lambda_2(\omega)) = 
\sum_{k=0}^{N_{\xi}-1} \lambda^{(2)}_k \mathbb{E}(\psi_k(\xi(\omega)))= \lambda^{(2)}_0,
 \]
due to the orthogonality of the basis polynomials $\{\psi_k\}.$ We remark here also that this mean value
of the second eigenvalue is quite close to the eigenvalue computed from the associated deterministic 
problem, i.e., $1.4166$. If we increased the order of the Legendre polynomials, then the INBM would 
converge to the exact deterministic value. However, this would come at a higher computational expense, as the 
quantity $N_{\xi}$ would also need to be increased accordingly. This kind of observation has earlier been extensively verified by the authors
 in the context of linear systems arising from PDEs with uncertain inputs \cite{BOS14, AO16}.
 
\begin{center}
\captionof{table}{{The first eight coefficients  of the  spectral expansion gPCE of the second eigenvalue with  using INBM.
Here, $k$ stands for the index of the basis function in the expansion (\ref{lambdas}).}}
\vspace{-4.5mm} 
\begin{tabular}{l|llllllllll}
\hline
$k$& $ 0$ &$1$&  $2$ & $3$&  $ 4$ &$5$& $6$&  $7$ \\
\hline
\hline
$\lambda^{(2)}_k$ & $1.4121$ & $0.5492$ &  $0.7009$ & $0.5492$&  $ -0.02013$  &$0.0952$& $-0.03537$ & $ 0.0594$\\
\hline
\end{tabular}
\label{tractab}
\end{center}
\vspace{-4.5mm} 

Next,  in Tables \ref{bicgstab_eta} and \ref{bicgstab_eta2}, we show the performance of the INBM solver 
in the two cases; that is, for dim$(\mathcal{J})=10,500$ and dim$(\mathcal{J})=82,468,050$.
Here, the Newton equations (cf. (\ref{NEq}))
are solved using low-rank BiCGstab, as well as using the standard preconditioned BiCGstab  method which we have denoted as full model (FM),
 that is, without low-rank truncation.
The CPU times reported are  the total time it took the solver to compute the spectral
coefficients of the eigenpair $(\lambda_2(\omega), \varphi_2(\omega)).$ Here, for each choice of the 
forcing terms $\{\eta_k\}$ discussed in Section \ref{inm}, we report inexact Newton steps (INS), backtracks per inexact Newton step (BINS), 
total BiCGstab iterations (iter),  total CPU times (t) in seconds, ranks of the solution (R), memory in kilobytes of the low-rank solution (LR) 
and full method solution (FM). 
By the memory requirement of a low-rank solution $X=WV^T,$ we mean the sum of the two separate computer memories occupied by its factors $W$ and $V^T,$
since $X$ is computed and stored in this format, unlike the solution from FM.
From the two tables, we see that for this problem, the performance of the INBM algorithm is independent of the 
the choice of the forcing terms $\{\eta_k\}.$ In particular, Table \ref{bicgstab_eta} shows that the algorithm could compute the solution within a
few seconds in the first case. Besides, the INBM algorithm reduces the storage requirement of the solution to one-quarter of 
memory required to compute the full solution. In fact, as shown in \cite{BOS14, BOS15a}, for a fixed $N_{\xi},$
low-rank Krylov solvers typically provide more storage benefits as $N_x \rightarrow \infty.$

Finally, as in the first case, we have also computed only the second stochastic eigenvalue  $\lambda_2(\omega) $ of the matrix $\mathcal{A}(\omega)$ 
(cf. (\ref{rmatrix})) for the case where dim$(\mathcal{J})=82,468,050$. Again, the mean $\lambda^{(2)}_0$ of this stochastic eigenvalue corresponds to the 
second eigenvalue of the deteministic matrix $A_0$, which in this case is $0.003$.
Note in particular from Table~\ref{bicgstab_eta2} that with the FM, MATLAB  indeed fails as the 
size of the Jacobian matrix $\mathcal{J}$ at each inexact Newton step is now increased to more than $82$ million degrees of freedom. 
Yet, INBM handles this task in about 200 minutes;  that is, roughly $6$ minutes per Newton step.
Here, the solution from FM terminates with `out of memory', which we have denoted as `OoM'.

\vspace{.3mm}
\begin{center}
\captionof{table}{Performance of the INBM solver for dim$(\mathcal{J})=10,500$ with $\sigma_a=0.01$. Here,
$I, II,$ and $III$ represent the  different forcing  parameter  choices (\ref{eta1}), (\ref{eta2}), and $(b)$ in Section \ref{inm}.
}
\begin{tabular}{l|llllllll}
\hline
$\eta_k$ & INS & BINS &$\#$ iter  & t & R& mem(LR) & mem(FM)\\
\hline
\hline
$I$& 22  & 1.5        & 22  &    16.5   & 9    & 18.7  & 84\\
$II$&22  & 1.5        & 22  &    17.5   & 10   & 20.8  & 84\\
$III$&22  & 1.5       & 23  &    17.2   & 10   & 20.8  & 84\\
\hline
\end{tabular}
\label{bicgstab_eta}
\end{center}

\vspace{.3mm}
\begin{center}
\captionof{table}{Performance of the INBM solver for dim$(\mathcal{J})=82,468,050$ with $\sigma_a=0.1$. Here,
$I, II,$ and $III$ represent the  different forcing  parameter  choices (\ref{eta1}), (\ref{eta2}), and $(b)$ in Section \ref{inm}.
}
\begin{tabular}{l|lllllll}
\hline
$\eta_k$ & INS & BINS &$\#$ iter  & t   & R& mem(LR) & mem(FM) \\
\hline
\hline
$I$   & 34    & 3.6        & 39  &      12123.4 & 51 & 156551.7 & OoM\\
$II$  & 32    & 3.5        & 43  &      12112.8 & 51 & 156551.7 & OoM\\
$III$ & 33    & 3.5        & 42  &      12200.1 & 51 & 156551.7 & OoM\\
\hline
\end{tabular}
\label{bicgstab_eta2}
\end{center}

\section{Conclusions}
\label{conc}
In computational science and engineering, there are certain problems of growing interest for which random matrices are considered
as random perturbations of finite-dimensional operators. These random matrices are usually not obtained from a
finite-dimensional representation of a partial differential operator, and in a number of interesting
cases, closed-form expressions of the statistical moments and probability density functions of their
eigenvalues and eigenvectors are available; see e.g., \cite{S05}. The matrices of interest in the present paper, on
the other hand, are the result of a finite-dimensional approximation of an underlying continuous
system and their stochasticity is intrinsically tied to the uncertainty in the parameters of this system. For
such systems, closed-form expressions are generally not available for the solution of the SEVPs.

In this paper, we have presented a low-rank Newton-type algorithm for approximating the eigenpairs of  SEVPs.
The numerical experiments confirm that the proposed solver can mitigate the 
computational complexity associated with solving high dimensional Newton systems in the considered SEVPs.
 More specifically, the low-rank approach guarantees significant storage savings
 \cite{BOS14, BOS15a, StoB13}, thereby enabling the 
 solution of  large-scale SEVPs that would otherwise be intractable.
 
%


\section*{Acknowledgement}
The authors would like to thank Sergey Dolgov  for fruitful discussions in the 
course of writing this paper. The work was performed while Martin Stoll was at the Max
Planck Institute for Dynamics of Complex Technical Systems.

\bibliographystyle{siam}
\bibliography{akwumref}

\end{document}